\pdfoutput=1
\documentclass[sigconf]{acmart}

\copyrightyear{2026}
\acmYear{2026}
\setcopyright{cc}
\setcctype{by}
\acmConference[ICAIF '26]{7th ACM International Conference on AI in Finance}{November 14--17, 2026}{Milan, AA, Italy}
\acmBooktitle{7th ACM International Conference on AI in Finance (ICAIF '26), November 14--17, 2026, Milan, AA, Italy}
\acmDOI{10.1145/3858063.3858709}
\acmISBN{979-8-4007-3087-0/2026/11}

\allowdisplaybreaks
\usepackage{nicefrac}
\usepackage{bm}
\usepackage[ruled,vlined]{algorithm2e}
\usepackage{enumitem}
\usepackage{xspace}

\setlist{topsep=2pt, partopsep=0pt, itemsep=2pt, parsep=0pt, leftmargin=*}
\setlist[enumerate]{label={\rm(\roman*)}}

\providecommand{\shuffle}{\mathbin{\sqcup\mkern-3.5mu\sqcup}}
\providecommand{\thicksim}{\sim}

\theoremstyle{remark}
\newtheorem{remark}{Remark}

\numberwithin{equation}{section}

\title{Signature Methods for Optimal Market Making}

\author{Alberto Gennaro}
\affiliation{\institution{UC Berkeley, IEOR}\city{Berkeley}\country{USA}}
\affiliation{\institution{MLCoE, JPMC}\city{New York}\country{USA}}
\email{alberto.gennaro@berkeley.edu}
\author{Thibaut Mastrolia}
\affiliation{\institution{UC Berkeley, IEOR}\city{Berkeley}\country{USA}}
\email{mastrolia@berkeley.edu}
\author{Francesca Primavera}
\affiliation{\institution{\'Ecole Polytechnique, CMAP\\ Engie Global Markets}\city{Paris}\country{France}}
\email{francesca.primavera@polytechnique.edu}

\begin{document}

\begin{abstract}
We propose a signature-based method to solve the optimal market-making problem under a mean-variance criterion. By exploiting signature linearization techniques, we reduce the market-making problem to a pseudo-linear optimization over the expected signature of an augmented market path, and we develop a signature algorithm named Sig-REINFORCE to learn the optimal bid and ask quotes. We test our method in two scenarios, in which market-order arrivals follow either a Poisson or a self-exciting Hawkes process, and we benchmark it against a Proximal Policy Optimization (PPO) baseline.
\end{abstract}

\ccsdesc[500]{Computing methodologies~Reinforcement learning}
\ccsdesc[500]{Mathematics of computing~Stochastic processes}
\ccsdesc[500]{Theory of computation~Continuous optimization}
\ccsdesc[500]{Applied computing~Economics}

\keywords{market making, signature methods, Hawkes processes, reinforcement learning, REINFORCE, policy gradient}

\maketitle

\section{Introduction}\label{sec:intro}

\subsection{The market making problem}
The seminal work of Ho and Stoll \cite{ho1981optimal} set the foundations of the market making problem, in which a trader optimally sets their bid/ask quotes on a limit order book, facing a continuous-time order flow and return risk on their inventory. It was then mathematically formalized in \cite{avellaneda2008high,stoikov2009option,Dealingrisk,euch2021optimal, ventre}, where the price of a risky asset follows a Bachelier model with Poisson order flow, and extended to option pricing in the Heston \cite{baldacci2021algorithmic} and rough Heston \cite{rosenbaum2022multi} models, multiple market makers \cite{baldacci2021optimal}, noise components \cite{barucci2025market}, and Hawkes processes \cite{jusselin2021optimal}; see also the monographs \cite{o1998market,cartea2015algorithmic}. Most of these studies cannot, however, reflect the complexity of both order flows and the fundamental price, nor cover all the empirical facts observed in financial markets: as Cont notes in \cite[Conclusion]{cont2001empirical}, existing models generally fail to reproduce all the statistical features of returns at once.

Model-free reinforcement learning methods for market making have been well developed in the last decades, from the foundations in \cite{nevmyvaka2006reinforcement} to the mean-variance criterion \cite{wang2020continuous}, multi-agent settings \cite{wang2025multi}, adversarial reinforcement learning \cite{zelman2024adversarial}, and adaptive market making anticipating closing auctions \cite{graf2026learning}; see \cite{beysolow2019market,hambly2023recent} for surveys. The two approaches have complementary drawbacks: model-based methods are rigorous and numerically tractable but inherently restrictive, since the perfect model does not exist, whereas RL methods absorb the specificity of any market but converge with difficulty in practice \cite{dulac2021challenges}.

We seek a middle ground: a mathematically rigorous method that retains the properties of the process rather than imposing a model. A natural tool for this is the signature of a stochastic process.

\subsection{Signature methods: background and motivation}
The signature of a path---the collection of all its iterated integrals---originates in \cite{CHEN:57} and plays a fundamental role in rough path theory \cite{L:98}. More recently, signatures have become a powerful tool in machine learning and mathematical finance for extracting informative features from time-series data, with applications to, e.g., optimal control \cite{BHRS:21,Bank2025StochasticCW, sigformer,aqsha2026solving}, data generation \cite{Buhler,NSSBWL:21,liao2024sig}, and financial modeling \cite{AbiJaber_Fouriersig,Eduardopathdependendsig,ASS:20,Alvaro_Immanol:22,CGS:22,KLP:20,LNP:20, pricingBlanka}.

From an approximation-theoretic perspective, iterated integrals play the role of monomials on path space. Three properties matter here: products of signature-linear functionals remain signature-linear via the shuffle product \eqref{eq:shuffle}; linear functionals of the signature are dense in suitable spaces of continuous path functionals \cite{KLP:20,CPS:22}; and sufficiently regular non-anticipative functionals admit expansions in signature coordinates \cite{CGP:24,DT:23}.
These properties make signature-linear functionals both expressive and tractable as policy classes for stochastic control. Since only a linear readout must be optimized, a broad class of non-Markovian control problems can be reduced to finite-dimensional optimization. Signature methods have consequently been applied successfully to a range of problems, including portfolio optimization and algorithmic trading (see e.g.,~\cite{Matteo_randomportfolio, Alvaro_Immanol:22, cuchiero2025signature,blankamarkovixsig,KLP:20,LNP:20}).

While signature methods have become an established tool for continuous-path systems, their extension to c\`adl\`ag paths is comparatively recent, and practical applications remain largely unexplored. This work develops a signature-based approximation framework for path-dependent market making and, to the best of our knowledge, provides the first applied use of signatures for stochastic control problems involving jump processes. The market-making setting considered here offers a natural and important testbed for this emerging theory, since order-flow dynamics are inherently discontinuous and are typically modeled by counting processes.

The signature of a c\`adl\`ag\xspace rough path was introduced in \cite{FS:17}, where the classical iterated-integral construction is modified to account for the interaction between continuous and jump components. The corresponding universal approximation theorem for c\`adl\`ag\xspace paths was subsequently established in \cite{CPS:22}. This result is essential for the present work, as the order-flow processes driving market-making dynamics are inherently discontinuous.
Specifically, the market maker posts ask and bid offsets $\delta^a$ and $\delta^b$ in response to buy and sell market orders represented by counting processes $N^a$ and $N^b$.

Our approach is based on a signature representation of path-dependent controls. More precisely, we parameterize the ask/bid controls as
\begin{center}
    $ \delta_t^a = \langle \bm{\alpha}^a, \mathbb{Y}_t \rangle,
\qquad
\delta_t^b = \langle \bm{\alpha}^b, \mathbb{Y}_t \rangle$
\end{center}
where $\bm{\alpha}^a$ and $\bm{\alpha}^b$ are trainable tensors supported on words up to a prescribed truncation level and  $\mathbb{Y}_t$ denote the signature of an augmented market path at time $t$, incorporating time, price information, quadratic-variation covariates, and order-flow history (see Section \ref{sec:signaturemethods}). This choice is not merely a computational convenience. By the c\`adl\`ag approximation theorem of \cite{CPS:22}, signature-linear functionals form a universal class for sufficiently regular non-anticipative controls, implying that any such control can be approximated arbitrarily well within this parameterization.

A second advantage of the signature framework is that it renders the optimization objective highly tractable. Using the algebraic structure of signatures, stochastic integrals of linear functionals of the signature against counting processes can be represented, after an explicit \emph{tilde transform} (Definition \ref{def:tildetrans}), as linear functionals of the terminal signature. Nonlinear interactions are handled through the shuffle product identity, which expresses products of signature features as elements of the same algebra. As a result, the realized trading payoff admits a finite-dimensional polynomial representation in the control parameters.

Following this reformulation, the principal source of complexity is no longer the payoff functional itself, but the dependence of the controlled probability law on the policy through the order-arrival intensities. We address this challenge using a likelihood-ratio gradient estimator, which enables direct optimization over a rich class of path-dependent controls while preserving computational efficiency.

\subsection{Contributions and outline}
The contributions of this work are fourfold.  First, we formulate the first signature-linear class of path-dependent market-making controls for semimartingale price and point-process order-flow dynamics. As far as we know, this is also the first application of signature methods to optimal control with jump processes.

Second, we provide a comprehensive and explicit transformation which rewrites the relevant stochastic integrals driven by continuous and counting processes (modeling the price and order-arrival components respectively) as linear functionals of the signature of an augmented path.

Third, we show how the mean--variance objective investigated can be reduced, via the shuffle product, to an algebraic functional of the expected terminal signature under the controlled law.

Fourth, we combine this representation with a likelihood-ratio gradient for spread-dependent point-process intensities and score functions.

The outcome is a fully tractable Sig-REINFORCE algorithm taking as input the signature of the observed data and returning the bid/ask quotes maximizing the mean-variance objective.

Section \ref{sec:preliminaries} recalls the notation and main results from signature theory and presents the general problem together with the linearization results. Section \ref{sec:mv_objective} applies them to the mean-variance problem, combining shuffle reduction with the likelihood-ratio gradient, and provides the Sig-REINFORCE algorithm. Section \ref{sec:numerics} reports numerical results in Poisson and Hawkes order-flow environments, benchmarking Sig-REINFORCE against the on-policy PPO \cite{schulman2017proximal}. 

\section{Signature market making optimization}\label{sec:preliminaries}

\subsection{Notations}\label{sec:notationsig}
Throughout, we consider a probability space $(\Omega,\mathcal F,\mathbb P)$ such that $\Omega=\mathcal C([0,T];\mathbb R)\times \mathbb D([0,T];\mathbb N)\times\mathbb D([0,T];\mathbb N)$ where $\mathcal C([0,T];\mathbb R)$ and $\mathbb D([0,T];\mathbb N)$ denote the space of continuous and c\`adl\`ag functions from $[0,T]$ into $\mathbb R$ or $\mathbb N$ respectively.
 Fix $d\in\mathbb N$. The extended tensor algebra over $\mathbb{R}^d$ (the Euclidean space) is defined by
\[
T((\mathbb{R}^d))
:=
\left\{
\mathbf{u}=(\mathbf{u}^{(0)},\mathbf{u}^{(1)},\dots,\mathbf{u}^{(n)},\dots)
\ \middle|\ 
\mathbf{u}^{(n)}\in (\mathbb{R}^d)^{\otimes n}
\right\},
\]
where $(\mathbb{R}^d)^{\otimes n}$ denotes the $n$-fold tensor product of $\mathbb{R}^d$, with the convention $(\mathbb{R}^d)^{\otimes 0}:=\mathbb{R}$. We equip
$T((\mathbb{R}^d))$ with the standard addition $+$, tensor multiplication $\otimes$, and scalar multiplication. For $K\in\mathbb{N}$, the truncated tensor algebra is defined by
\[
T^{K}(\mathbb{R}^d)
:=
\left\{
\mathbf{u}=(\mathbf{u}^{(0)},\mathbf{u}^{(1)},\dots,\mathbf{u}^{(K)})
\ \middle|\ 
\mathbf{u}^{(n)}\in (\mathbb{R}^d)^{\otimes n}\ \text{for }n\leq K
\right\}.
\]
Let $\pi_n:T((\mathbb{R}^d))\to(\mathbb{R}^d)^{\otimes n}$ and $\pi_{\leq K}:T((\mathbb{R}^d))\to T^{K}(\mathbb{R}^d)$ be such that 
\[
\pi_n(\mathbf{u})=\mathbf{u}^{(n)}, \quad \pi_{\leq K}(\mathbf{u})=\mathbf{u}^{K}:=(\mathbf{u}^{(n)})_{n=0}^{K},
\]
for $\mathbf{u}\in T((\mathbb{R}^d))$.
Let $J=(j_1,\dots,j_n)$ be a multi-index with entries in $\{1,\dots,d\}$ and $\epsilon_1,\ldots,\epsilon_d$ the canonical basis of $\mathbb R^d$. We set $|J|:=n$, 
\[
\epsilon_J:=\epsilon_{j_1}\otimes \epsilon_{j_2}\otimes\dots\otimes \epsilon_{j_n}.
\]
Furthermore, we denote by $\epsilon_\emptyset$ the basis element of $(\mathbb{R}^d)^{\otimes 0}$ and set $|\emptyset|:=0$. We also set $J':=(j_1,\dots,j_{n-1})$ for $n>1$, $J':=\emptyset$ for $n=1$, and $J'':=(J')'$ for $n>1$, with the convention $\epsilon_{J''}=0$ for $n=1$. For two multi-indices $J\in\{1,\dots,d\}^{|J|}$ and $L\in\{1,\dots,d\}^{|L|}$, and two indices of length $1$, $a,b\in\{1,\dots,d\}$, the shuffle product $\shuffle$ is defined recursively by
\begin{align}\label{eq:shuffle}
    &J\shuffle \emptyset=\emptyset\shuffle J=J,\\
    &(J,a)\shuffle (L,b)=((J\shuffle (L,b)),a)+(((J,a)\shuffle L),b),\nonumber
\end{align}
where $(J,a)$ denotes the concatenation of multi-indices. Given $\mathbf{x}\in T((\mathbb{R}^d))$, we write $\mathbf{x}_J:=\langle \mathbf{x},\epsilon_J\rangle$, and for each $\mathbf{u}\in T^K(\mathbb{R}^d)$, we set
\begin{equation}\label{eq:notation_pairing}
    \langle \mathbf{u},\mathbf{x}\rangle:=\sum_{|J|=0}^K\mathbf{u}_J \mathbf{x}_J \in \mathbb{R}.
\end{equation}
Finally, we define
\begin{equation}\label{eq:G_group}
\begin{aligned}
    G((\mathbb{R}^d))
    :=
    \bigl\{
    &\mathbf{x}\in T((\mathbb{R}^d))
    \colon
    \langle \epsilon_\emptyset,\mathbf{x}\rangle=1,\\
    &\langle \epsilon_J\shuffle \epsilon_L,\mathbf{x}\rangle
    =
    \langle \epsilon_J,\mathbf{x}\rangle \langle \epsilon_L,\mathbf{x}\rangle,
    \ \text{for all }J,L
    \bigr\}.
\end{aligned}
\end{equation}

\subsection{Signature of c\`adl\`ag semimartingales}
In this section, we recall the notion of the signature of c\`adl\`ag semimartingale, as introduced in \cite{FS:17}.
\begin{definition}
Fix $T>0$. For a semimartingale $X=(X_t)_{t\in[0,T]}$, we call the $G((\mathbb{R}^d))$-valued process $\mathbb{X}=(\mathbb{X}_t)_{t\in[0,T]}$ the \emph{signature of $X$} if for all $t\in[0,T]$, $\langle \epsilon_\emptyset,\mathbb{X}_t\rangle=1$ and
\begin{align*}
    \langle \epsilon_J,\mathbb{X}_t\rangle
    ={}&
    \int_0^t \langle \epsilon_{J'},\mathbb{X}_{s^-}\rangle\, dX_s^{j_{|J|}}\,\mathbf{1}_{\{|J|>0\}}
    \\
    &
    +\frac{1}{2}\int_0^t \langle \epsilon_{J''},\mathbb{X}_{s^-}\rangle\, d[X^{j_{|J|-1}},X^{j_{|J|}}]^c_s\,\mathbf{1}_{\{|J|>1\}}
    \\
    &
    +\sum_{0<s\leq t}\sum_{\epsilon_{J_1}\otimes \epsilon_{J_2}=\epsilon_J}
    \frac{1}{|J_2|!}\,\mathbf{1}_{\{|J_2|>1\}}
    \langle \epsilon_{J_1},\mathbb{X}_{s^-}\rangle
    \langle \epsilon_{J_2},(\Delta X_s)^{\otimes |J_2|}\rangle,
\end{align*}
for each $|J|>0$.
\end{definition}

\subsection{Market making problem}
We consider a market maker trading a financial instrument whose mid-price process is denoted
by $S$, over an intraday horizon $[0,T]$; interest rates are neglected throughout. At each
time $t$, the market maker posts an ask price $S_t^a$ and a bid price $S_t^b$, or
equivalently controls the half-spreads
\[
    \delta_t^a := S_t^a - S_t,
    \qquad
    \delta_t^b := S_t - S_t^b,
\]
which we refer to as the ask and bid offsets, respectively. Let $N^a$ and $N^b$ be counting processes representing the cumulative arrivals of market
orders hitting the ask and the bid, respectively, with intensities $\lambda^a$ and
$\lambda^b$. We assume that
\begin{align}\label{eq:intensities}
    \lambda^a_t = h^a(\delta^a_{t^-}),
    \qquad
    \lambda^b_t = h^b(\delta^b_{t^-}),
\end{align}
for some positive valued functions $h^a,\, h^b$. The inventory process 
\[
    Q_t := N_t^b - N_t^a,
\]
tracks the position accumulated by the market maker up to time $t$.
The profit and loss generated over $[0,T]$ decomposes as
\[
    \int_0^T \!\bigl(S_t^a\,dN_t^a - S_t^b\,dN_t^b\bigr)
    +
    \bigl(Q_T S_T - Q_0 S_0\bigr)
    =
    Z_T^\delta + I_T,
\]
where, for $\delta := (\delta^a, \delta^b)$, $Z^\delta$ is the execution revenue, capturing the spread income earned on each
trade, and $I$ is the inventory P\&L, driven by fluctuations in the mid-price:
\[
    dZ_t^\delta := \delta_t^a\,dN_t^a + \delta_t^b\,dN_t^b,
    \qquad
    dI_t := Q_t\,dS_t.
\]
The market maker's problem is to choose the quiting strategy $\delta$
so as to maximize the expected utility of the total P\&L:
\begin{equation}\label{eq:optimalMM}
    V_0^{MM}
    =
    \sup_{\delta\,=\,(\delta^a,\,\delta^b)}
    \mathbb{E}\!\left[
        U\!\left(Z_T^\delta,\, I_T,\, (I_T)^2\right)
    \right],
\end{equation}
for some utility function $U$ to be specified.

\subsection{Signature methods}\label{sec:signaturemethods}
We show that, by restricting to the class of signature-linear offsets introduced 
below, the optimization problem~\eqref{eq:optimalMM} admits a tractable 
reformulation in terms of the signature of a state process, which 
summarizes the relevant market information up to the current time.

\subsubsection{State process and notation.}
Without loss of generality, we assume that $\Delta N_t^a\,\Delta N_t^b = 0$ for all $t$ and define the state process
\[
\begin{aligned}
Y_t
:=
\Bigl(
&t,\, S_t,\, [S]_t,
N_t^a,\,{\textstyle\sum_{s\leq t}}(\Delta N_s^a)^2,\,\dots,\,
{\textstyle\sum_{s\leq t}}(\Delta N_s^a)^M,\\
&N_t^b,\,{\textstyle\sum_{s\leq t}}(\Delta N_s^b)^2,\,\dots,\,
{\textstyle\sum_{s\leq t}}(\Delta N_s^b)^M
\Bigr),
\end{aligned}
\]
for some $M \geq 2$, whose components we denote by
\[
Y_t^0 := t, \quad
Y_t^s := S_t, \quad
Y_t^{qs} := [S]_t, \quad
Y_t^a := N_t^a, \quad
Y_t^b := N_t^b, 
\]
\[
Y_t^{a^m} := \sum_{s \leq t}(\Delta N_s^a)^m, \quad
Y_t^{b^m} := \sum_{s \leq t}(\Delta N_s^b)^m,
\]
for $2 \leq m \leq M$. Multi-indices have entries in the alphabet
\[
\mathcal{A} := \{0, s, qs, a, a^2, \dots, a^M, b, b^2, \dots, b^M\},
\]
with $|\mathcal{A}| = 3 + 2M$, using the notation of Section~\ref{sec:notationsig}.

\subsubsection{Linearization of stochastic integrals.}
\begin{definition}\label{def:tildetrans}
For each multi-index $J$ with entries in $\{0,s,a,b\}$ and each
$j \in \{0,s,a,b\}$, define
\begin{align*}
    (\epsilon_J;\epsilon_j)^{\thicksim}
    :={}&
    \epsilon_J \otimes \epsilon_j
    - \frac{1}{2}\,\epsilon_{J'} \otimes \epsilon_{qs}\,\mathbf{1}_{\{i_{|J|}=j=s\}}
    \\
    &+ \sum_{\epsilon_{J_1} \otimes \epsilon_{J_2} = \epsilon_J}
    \bm{c}(|J_2|)\,(\epsilon_{J_1} \otimes \epsilon_{a^{|J_2|+1}})\,
    \mathbf{1}_{\{j=a\}}\,\mathbf{1}_{\{J_2=(a,\dots,a)\}}
    \\
    &+ \sum_{\epsilon_{J_1} \otimes \epsilon_{J_2} = \epsilon_J}
    \bm{c}(|J_2|)\,(\epsilon_{J_1} \otimes \epsilon_{b^{|J_2|+1}})\,
    \mathbf{1}_{\{j=b\}}\,\mathbf{1}_{\{J_2=(b,\dots,b)\}},
\end{align*}
for $\bm{c}(r) := \sum_{k=1}^r (-1)^k \bm{c}(r,k)$, 
$\bm{c}(r,k) := \sum_{J \in \mathcal{I}(r,k)} \prod_{i=1}^k \frac{1}{(j_i+1)!}$,
with $\mathcal{I}(r,k)$ denoting the set of multi-indices $J \in \mathbb{N}^k$,
with $j_i \geq 1$ for every $i$, summing to $r$.
\end{definition}
\begin{remark}
Under the additional assumption $[S]_t = \sigma^2 t$, the term
$\epsilon_{J'} \otimes \epsilon_{qs}$ collapses to $\sigma^2\,\epsilon_{J'} \otimes \epsilon_0$.
If furthermore all jumps are of unit size, $\Delta N_t^a = \Delta N_t^b = 1$
at every jump time, then $\epsilon_{a^{|J_2|+1}}$ and $\epsilon_{b^{|J_2|+1}}$ reduce
to $\epsilon_a$ and $\epsilon_b$, yielding a considerably simpler expression.
Both simplifications are exploited in the subsequent sections.
\end{remark}
The following proposition is the cornerstone of our approach: it shows that
It\^o integrals of signature-linear functionals remain signature-linear.
The proof is a simple adaptation of Proposition~4.10 in~\cite{CPS:22}.
\begin{proposition}\label{prop:LinearIntegral}
Fix $n \in \mathbb{N}$, $J \in \{0,s,a,b\}^n$, and $j \in \{0,s,a,b\}$.
If $M \geq n+1$, then, for all $t\in [0,T]$,
\[
\int_0^t \langle \epsilon_J, \mathbb{Y}_{s^-} \rangle\, dY_s^j
=
\langle (\epsilon_J;\epsilon_j)^{\thicksim}, \mathbb{Y}_t \rangle.
\]
\end{proposition}
By linearity of the stochastic integral, this extends immediately to
$\mathbf{u} \in T^K(\mathbb{R}^4)$: for every $j \in \{0,s,a,b\}$,
\[
\int_0^t \langle \mathbf{u}, \mathbb{Y}_s \rangle\, dY_s^j
=
\langle (\mathbf{u};\epsilon_j)^{\thicksim}, \mathbb{Y}_t \rangle,
\qquad
(\mathbf{u};\epsilon_j)^{\thicksim} := \sum_{|J| \geq 0} u_J\,(\epsilon_J;\epsilon_j)^{\thicksim}.
\]

\subsubsection{Signature-linear offsets and problem reformulation.}
We introduce the class of signature-linear offsets, i.e., controls parametrized
as linear functionals of the signature of $Y$, and show how
Proposition~\ref{prop:LinearIntegral} yields a tractable reformulation
of~\eqref{eq:optimalMM}. The restriction to this class is natural and well-motivated:
by the universal approximation properties of the signature of c\`adl\`ag paths (see \cite{CPS:22}),
any sufficiently regular control can be approximated arbitrarily well by a
signature-linear one.
\begin{definition}[Lin-Sig-Offsets]\label{def:linsigoffsets}
Let $K \in \mathbb{N}$ be fixed. The space of ($K$-level) \emph{signature-linear offsets}
is $\mathcal{D}^{\mathrm{sig}} := \mathcal{D}^{\mathrm{sig},a} \times
\mathcal{D}^{\mathrm{sig},b}$, where
\begin{align*}
    \mathcal{D}^{\mathrm{sig},a}
    &:=
    \Bigl\{
        \delta_\cdot^a = \langle \bm{\alpha}^a, \mathbb{Y}_\cdot \rangle
        :=
        {\textstyle\sum_{|J| \leq K}} \bm{\alpha}_J^a \langle \epsilon_J, \mathbb{Y}_\cdot \rangle
        \;\Big|\;
        \bm{\alpha}_J^a = 0 \text{ if } J \cap \bar{\mathcal{A}} \neq \emptyset
    \Bigr\},
    \\
    \mathcal{D}^{\mathrm{sig},b}
    &:=
    \Bigl\{
        \delta_\cdot^b = \langle \bm{\alpha}^b, \mathbb{Y}_\cdot \rangle
        :=
        {\textstyle\sum_{|J| \leq K}} \bm{\alpha}_J^b \langle \epsilon_J, \mathbb{Y}_\cdot  \rangle
        \;\Big|\;
        \bm{\alpha}_J^b = 0 \text{ if } J \cap \bar{\mathcal{A}} \neq \emptyset
    \Bigr\},
\end{align*}
where $\bar{\mathcal{A}}:=\mathcal{A} \setminus \{0,s,a,b\}$.
\end{definition}
Now we show that restricting to Lin-Sig-Offsets allows us to recast the optimization problem~\eqref{eq:optimalMM} as a tractable optimization problem.
First, note that as a consequence of Proposition~\ref{prop:LinearIntegral}, for any
 $\delta \in \mathcal{D}^{\mathrm{sig}}$, both $Z^\delta$ and
$I$ are linear functionals of $\mathbb{Y}$.
\begin{corollary}\label{cor:ZIlinear}
For any $\delta = (\langle \bm{\alpha}^a, \mathbb{Y}_\cdot \rangle,
\langle \bm{\alpha}^b, \mathbb{Y}_\cdot \rangle) \in \mathcal{D}^{\mathrm{sig}}$,
\[
Z_\cdot^\delta
=
\langle (\bm{\alpha}^a;\epsilon_a)^{\thicksim} + (\bm{\alpha}^b;\epsilon_b)^{\thicksim},\mathbb{Y}_\cdot \rangle, \quad I_\cdot
=
\langle (\epsilon_b - \epsilon_a;\epsilon_s)^{\thicksim}, \mathbb{Y}_\cdot \rangle.
\]
\end{corollary}
\begin{remark}\label{remark:2sources}
Two sources of dependence on $\delta$ are present. The first is explicit: the coefficients of 
$(\bm{\alpha}^a, \bm{\alpha}^b)\in T^K(\mathbb{R}^4)\times T^K(\mathbb{R}^4) $ appear directly in $Z^\delta$.
The second is implicit: through the intensity
condition~\eqref{eq:intensities}, the control $\delta$ also governs the law of $Y$,
hence of its signature $\mathbb{Y}$, which enters both $Z^\delta$ and
$I$. 
\end{remark}
To conclude, Corollary~\ref{cor:ZIlinear} allows us to recast~\eqref{eq:optimalMM} as
\begin{equation}\label{eq:optimalMMsig}
    \begin{aligned}
    &\sup_{\bm{\alpha}^a,\bm{\alpha}^b}
    \mathbb{E}\Big[
    \tilde{U}\Big(
    \langle (\bm{\alpha}^a;\epsilon_a)^{\thicksim}+(\bm{\alpha}^b;\epsilon_b)^{\thicksim},\mathbb{Y}_T \rangle, \langle (\epsilon_b-\epsilon_a;\epsilon_s)^{\thicksim},\mathbb{Y}_T \rangle
    \Big)
    \Big],
\end{aligned}
\end{equation}
where $\tilde{U}(\zeta,\iota) = U(\zeta,\iota,\iota^2)$. If $U$ is polynomial, or well approximated by one, by the shuffle property
of the signature~\eqref{eq:G_group}, the objective is then polynomial in the coefficients of 
$(\bm{\alpha}^a, \bm{\alpha}^b) $, with all dependence on the control law
encoded in
$\mathbb{E}[\mathbb{Y}_T]$.

\section{Mean-variance objective}\label{sec:mv_objective}
We specialize the general framework to a mean-variance objective and study the resulting problem via signature methods. This setting underlies the experiments of Section~\ref{sec:numerics}.
For any quoting strategy $\delta$, the objective function takes the form
\begin{equation}\label{eq:mv_objective}
    J(\delta)
    :=
    \mathbb{E}\!\left[Z_T^\delta\right]
    -\eta\,\mathbb{E}\!\left[\bigl(I_T\bigr)^2\right],
    \qquad \eta > 0,
\end{equation}
where, recall that
\[
dZ_t^\delta = \delta_t^a\,dN_t^a + \delta_t^b\,dN_t^b,
\qquad
dI_t = Q_t\,dS_t,
\qquad
Q_t = N_t^b - N_t^a.
\]

\begin{remark}
We call this criterion mean-variance referring to maximizing the mean of the profit and loss generated by the market making trading activities ($Z^\delta$) while minimizing the variance of the final  inventory $I_T$. Since $I$ is a martingale, it is equivalent to a quadratic penalty weighted by the risk aversion parameter $\eta>0$.  
\end{remark}
\subsection{Problem reformulation via signature methods}\label{sec:mv_signature_form}
To lighten notation, we set $\bm{\alpha}:=(\bm{\alpha}^a,\bm{\alpha}^b)$ and
\begin{equation*}
    \ell_Z(\bm{\alpha})
    :=
    (\bm{\alpha}^a;\epsilon_a)^{\thicksim}+(\bm{\alpha}^b;\epsilon_b)^{\thicksim},
    \qquad
    \ell_I
    :=
    (\epsilon_b-\epsilon_a;\epsilon_s)^{\thicksim}.
\end{equation*}
By Corollary~\ref{cor:ZIlinear} and the shuffle identity~\eqref{eq:G_group}, for the signature-linear quotes $\delta^{a/b}_\cdot=\langle\bm{\alpha}^{a/b},\mathbb{Y}_\cdot\rangle$ the objective~\eqref{eq:mv_objective} becomes
\begin{equation}\label{eq:mv_objective_signature_inner_product}
    J(\bm{\alpha})
    =
    \Bigl\langle
    L(\bm{\alpha}),\,
    \mathbb{E}_{\bm{\alpha}}[\mathbb{Y}_T]
    \Bigr\rangle,
    \quad
    L(\bm{\alpha})
    :=
    \ell_Z(\bm{\alpha})-\eta\,(\ell_I \shuffle \ell_I).
\end{equation}

Equation~\eqref{eq:mv_objective_signature_inner_product} reveals a \emph{pseudo-linear structure}: $J(\bm{\alpha})$ is linear in the expected signature, but $\bm{\alpha}$ enters through two terms (Remark~\ref{remark:2sources}):
\begin{enumerate}
    \item the functional $L(\bm{\alpha})$, which is explicit in $\bm{\alpha}$;
    \item the expected signature, which depends implicitly on $\bm{\alpha}$ since the Lin-Sig-Offsets drive the arrival intensities via~\eqref{eq:intensities}, and hence the law of the whole trajectory, which is why we write $\mathbb{E}_{\bm{\alpha}}[\mathbb{Y}_T]$.
\end{enumerate}
This implicit dependence does not obstruct gradient-based optimization: the gradient of $\mathbb{E}_{\bm{\alpha}}[\mathbb{Y}_T]$ with respect to $\bm{\alpha}$ is exactly what the REINFORCE estimator provides (Algorithm~\ref{alg:reinforce_sigmm}).

\subsection{Likelihood-ratio gradient}\label{sec:LR_gradient}
We rewrite the objective \eqref{eq:mv_objective_signature_inner_product} as
\begin{align}\label{eq:J}
    J(\bm{\alpha}) = \mathbb{E}_{\bm{\alpha}}\!\bigl[R(\omega;\bm{\alpha})\bigr],
\end{align}
where $R(\omega;\bm{\alpha})
    :=
    \bigl\langle \ell_Z(\bm{\alpha}),\mathbb{Y}_T(\omega)\bigr\rangle
    -\eta\bigl\langle \ell_I \shuffle \ell_I,\mathbb{Y}_T(\omega)\bigr\rangle$ is the pathwise reward. As highlighted in Section~\ref{sec:mv_signature_form}, $J$ depends on $\bm{\alpha}$ through two distinct channels: directly, via the reward $R(\omega;\bm{\alpha})$, and implicitly, via the controlled law $\mathbb{P}_{\bm{\alpha}}$, which redistributes probability mass across trajectories as the offsets vary. The latter dependence cannot be handled by a pathwise derivative, since the parameter shifts the underlying measure rather than a single realization; this is what motivates the likelihood-ratio approach described hereafter.

Suppose the family $(\mathbb{P}_{\bm{\alpha}})_{\bm{\alpha}}$ is dominated by a common reference measure $\mathbb{P}^\circ$ on path space, with Radon--Nikodym density $p_{\bm{\alpha}}(\omega) := d\mathbb{P}_{\bm{\alpha}}/d\mathbb{P}^\circ$ that is strictly positive and differentiable in $\bm{\alpha}$ for $\mathbb{P}^\circ$-a.e.\ $\omega$, and that both $R(\omega;\bm{\alpha})\,p_{\bm{\alpha}}(\omega)$ and $R(\omega;\bm{\alpha})\,\nabla_{\bm{\alpha}}p_{\bm{\alpha}}(\omega)$ admit a locally dominated $\mathbb{P}^\circ$-integrable envelope. Differentiating under the integral sign and applying the score-function identity $\nabla_{\bm{\alpha}}p_{\bm{\alpha}}(\omega) = p_{\bm{\alpha}}(\omega) 
\nabla_{\bm{\alpha}}\log p_{\bm{\alpha}}(\omega) $ yields
\begin{equation}\label{eq:grad_decomp}
    \nabla_{\bm{\alpha}} J(\bm{\alpha})
    =
    \underbrace{
        \mathbb{E}_{\bm{\alpha}}\!\left[
            \nabla_{\bm{\alpha}} R(\omega;\bm{\alpha})
        \right]
    }_{\text{explicit algebraic term}}
    +
    \underbrace{
        \mathbb{E}_{\bm{\alpha}}\!\left[
            R(\omega;\bm{\alpha})
            \nabla_{\bm{\alpha}}\log p_{\bm{\alpha}}(\omega)
        \right]
    }_{\text{likelihood-ratio term}}.
\end{equation}
The first term arises from the direct dependence of the execution revenue on $\bm{\alpha}$ through the tilde representation, and admits the closed-form expression
\[
    \nabla_{\bm{\alpha}} R(\omega;\bm{\alpha})
    =
    \bigl\langle \nabla_{\bm{\alpha}} \ell_Z(\bm{\alpha}),\mathbb{Y}_T(\omega)\bigr\rangle.
\]
The second accounts for the dependence of the controlled law $\mathbb{P}_{\bm{\alpha}}$ on $\bm{\alpha}$, and is estimated via the REINFORCE estimator detailed in the following section.

\subsection{Pseudo-algorithm}\label{sec:pseudo_algo}

We summarize a gradient-ascent scheme for maximizing \eqref{eq:J} on the coefficients of $\bm{\alpha}$ following the ideas in Section \ref{sec:LR_gradient}, which applies to all intensity models and price dynamics for which the score of the controlled law is available.

\begin{algorithm}[!ht]
\caption{Sig-REINFORCE}\label{alg:reinforce_sigmm}
\SetKwInput{KwIn}{Input}
\SetKwInput{KwReturn}{Return}
\SetKwFor{ForEach}{for each}{do}{end}
\SetAlgoLined
\DontPrintSemicolon
\SetNlSty{}{}{}
\SetAlCapFnt{\small}
\SetKwComment{Comment}{$\triangleright$\ }{}
\footnotesize
\KwIn{$K$, $B$, $(\rho_m)$, $\eta>0$, $G_{\max}$, $\varepsilon>0$, $m_{\max}$ (level, batch, rates, risk aversion, clip, tolerance, max iterations)}
\textbf{Init:} $\bm{\alpha}_0\in(T^K(\mathbb{R}^d))^2$ (e.g.\ $\mathbf{0}$), $m\leftarrow 0$\;
\Repeat{$\|\hat g_m\|<\varepsilon$ \emph{or} $m=m_{\max}$}{
  Draw $\omega^{[1]},\dots,\omega^{[B]}\sim\mathbb{P}_{\bm{\alpha}_m}$ \Comment*[r]{simulate}
  \ForEach{$i=1,\dots,B$}{
    $\mathbb{Y}_T^{[i]}\leftarrow$ signature of $\omega^{[i]}$ at $T$ \Comment*[r]{terminal sig.}
    $\phi^{[i]}_j\leftarrow\pi_{\leq K}(\mathbb{Y}^{[i]}_{\tau_j})\ \forall\tau_j$ \Comment*[r]{features at event times}
    $R^{[i]}\leftarrow\langle \ell_Z(\bm{\alpha}_m)-\eta\,\ell_I\shuffle\ell_I,\mathbb{Y}_T^{[i]}\rangle$ \Comment*[r]{reward}
    $g_{\mathrm{alg}}^{[i]}\leftarrow\langle\nabla_{\bm{\alpha}}\ell_Z(\bm{\alpha}_m),\mathbb{Y}_T^{[i]}\rangle$ \Comment*[r]{algebraic grad.}
    $s^{[i]}\leftarrow\nabla_{\bm{\alpha}}\log p_{\bm{\alpha}_m}(\omega^{[i]})$ \Comment*[r]{score via \eqref{eq:poisson_score_alpha} and $(\phi^{[i]}_j)_j$}
  }
  $\bar R_m\leftarrow\frac1B\sum_{i=1}^{B}R^{[i]}$ \Comment*[r]{baseline}
  $\hat g_m\leftarrow\frac1B\sum_{i=1}^{B}\bigl(g_{\mathrm{alg}}^{[i]}+(R^{[i]}-\bar R_m)\,s^{[i]}\bigr)$ \Comment*[r]{gradient estimate}
  $g\leftarrow\min\bigl(1,G_{\max}/\|\hat g_m\|\bigr)\,\hat g_m$ \Comment*[r]{clipping}
  $\bm{\alpha}_{m+1}\leftarrow\bm{\alpha}_m+\rho_m g$;\quad $m\leftarrow m+1$\;
}
\KwReturn{$\bm{\alpha}^\star=\bm{\alpha}_m$}
\end{algorithm}

We denote by $\bar R_m$ a baseline used for variance reduction. Note that, subtracting it in Algorithm~\ref{alg:reinforce_sigmm} does not bias the estimator, since the score $\nabla_{\bm{\alpha}}\log p_{\bm{\alpha}}$ has zero mean, and serves only to reduce its variance; see Section~\ref{sec:ppo}.

\section{Numerical experiments}\label{sec:numerics}
We finally turn to the implementation of the Sig-REINFORCE algorithm described in Section~\ref{sec:pseudo_algo}.
\subsection{Experimental design}\label{sec:numerics_setup}
We evaluate the method in two order flows regimes. For both, the underlying stock price follows a Bachelier model $S_t=S_0+\sigma W_t,$ with $S_0,\sigma>0$. For each training run we pre-sample $B$ independent Brownian paths $W^{[1]},\dots,W^{[B]}$ on a fine grid, one per trajectory of the batch, and freeze them.  The same $B$ paths are reused at every iteration, while the order flow is redrawn at every rollout. All out-of-sample comparisons of Sections~\ref{sec:Poisson}--\ref{sec:hawkes_numerics} use a second, independent set of price paths shared by Sig-REINFORCE, Sig-PPO and the expert, so that the methods differ only through their quotes (common random numbers). The consequences of this construction for the gradient estimator are discussed in Section~\ref{sec:ppo}.
The two order flows differ therefore only for the bid/ask arrivals:
\begin{enumerate}
    \item in the first, order flow is conditionally Poisson with spread-controlled intensities; this admits a closed-form mean--variance optimum \cite[Chapter 6]{cartea2015algorithmic} and \cite{Dealingrisk}, which we use as an \emph{expert} benchmark to check that Sig-REINFORCE recovers a known optimum (Section~\ref{sec:Poisson});
    \item in the second, order flow is a self-exciting Hawkes process (Section~\ref{sec:hawkes_numerics}): genuinely path-dependent and without a closed-form solution, yet still tractable, so that it stresses the path-dependent features of the signature policy.
\end{enumerate}  
In both regimes we benchmark against a PPO baseline (Section~\ref{sec:ppo}). Note that  the linearization of Section~\ref{sec:signaturemethods} is not tied to this choice of processes since the only requirement for Sig-REINFORCE is a counting process whose intensity is positive and differentiable in the quoted spreads, so that the score $\nabla_{\bm{\alpha}}\log p_{\bm{\alpha}}$ exists.

\begin{figure}[ht]
    \centering\includegraphics[width=0.9\columnwidth]{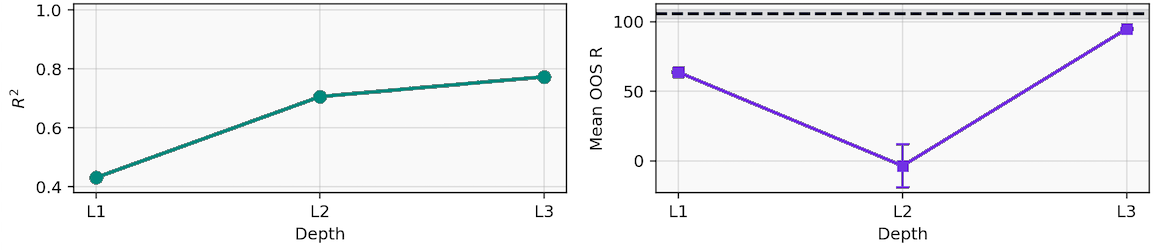}
    \caption{$R^2$ of the ridge fit to the expert policy and out-of-sample
    percentage reward of the fitted linear-in-signature rule, as a function of
    the signature truncation level $K$, for the Poisson order flow.}
    \label{fig:level_selection}
    \Description{Plot of the ridge-regression R-squared and the out-of-sample percentage reward as a function of the signature truncation level K, showing that K equal to 3 already recovers the expert reward.}
\end{figure}

We made few computational choices that we recall here.
\begin{enumerate}
   \item We restrict the policy alphabet to $\{0,a,b\}$. This is without loss of generality in this model: the price is a martingale and the optimal quotes do not depend on it, while $[S]_t=\sigma^2 t$ is deterministic, so the quadratic-variation channel is a multiple of the time channel. Moreover, since the quotes do not depend on the price, the order flow $(N^a,N^b)$ is independent of $W$, and the score $\nabla_{\bm{\alpha}}\log p_{\bm{\alpha}}$ is a functional of $(N^a,N^b)$ only. The price enters the objective only through $I_T$, and It\^o's isometry gives
$\mathbb{E}\big[I_T^2\,\big|\,N^a,N^b\big]=\sigma^2\!\int_0^T Q_t^2\,dt$, which does not depend on $W$. This is why the frozen-price batch construction of Section~\ref{sec:ppo} does not bias the gradient estimator.
    \item For the same reasons, we choose the truncation level $K=3$. A signature-linear policy is indexed by the words of length at most $K$ on $\{0,a,b\}$, whose number grows geometrically as $\sum_{j=0}^{K}3^{\,j}$ per side: each extra level enlarges the policy class, but also the dimension of the parameter vector. We select the smallest level that is expressive enough through a representability test that does not involve training. We generate paths under the closed-form expert \cite{cartea2015algorithmic}, compute the truncated signature of the augmented path, and ridge-regress the expert's quoted half-spread on the signature features. We report the in-sample $R^2$ of the fit and the out-of-sample percentage mean--variance reward of the resulting linear-in-signature rule (Figure~\ref{fig:level_selection}). The gain from $K=2$ to $K=3$ is substantial and $K=3$ already recovers the full expert reward, whereas $K=4$ adds parameters at virtually no benefit. This is consistent with the observations in \cite{hu2025deep}. We keep $K=3$ for Hawkes, where no closed-form expert is available.
\end{enumerate} 

\subsection{Sig-REINFORCE vs.\ Sig-PPO}\label{sec:ppo}
In this section we contrast the two learning procedures used in Sections~\ref{sec:Poisson}--\ref{sec:hawkes_numerics}. Both Sig-PPO and Sig-REINFORCE share the same environment, the same simulator and the same signature feature map. We deliberately pick PPO as the deep-RL baseline because it is an \emph{on-policy} method: each gradient update is computed from rollouts, the difference is only in how the gradient of the expected reward $ J(\bm{\alpha}) = \mathbb{E}_{\bm{\alpha}}\!\bigl[R(\omega;\bm{\alpha})\bigr]$ is estimated from rollouts.

\paragraph{Batch construction.}
Both methods are trained on the same frozen price paths $(W^{[i]})_{i\le B}$, while the order flow is redrawn at every rollout, so they optimize the sample-average objective $J_W(\bm{\alpha}):=\tfrac1B\sum_{i}\mathbb{E}_{\bm{\alpha}}\bigl[R\,\big|\,W^{[i]}\bigr]$, a Monte Carlo approximation of $J(\bm{\alpha})$ over the price paths. For Sig-REINFORCE, the $B$ trajectories of a batch are independent, so $\bar R_m$ and the estimator of Algorithm~\ref{alg:reinforce_sigmm} are averages of independent terms, and successive batches are i.i.d.\ conditionally on $(W^{[i]})_{i\le B}$.

\subsubsection{Sig-REINFORCE (path-wise score).}
We estimate $\nabla_{\bm{\alpha}} J(\bm{\alpha})$ via the explicit algebraic part and a single likelihood-ratio term applied to the whole trajectory of $(W,N^a,N^b)$ given by
\[
    \mathbb{E}_{\bm{\alpha}}\!\left[\bigl(R(\bm{\alpha})-b\bigr)\,
    \nabla_{\bm{\alpha}}\log p_{\bm{\alpha}}(W,N^a,N^b)\right]
    =\mathbb{E}\bigl[E_{\bm{\alpha}}(W)\bigr],
\]
where $E_{\bm{\alpha}}(\varsigma):=\mathbb{E}_{\bm{\alpha}}\bigl[(R(\bm{\alpha})-b)\,\nabla_{\bm{\alpha}}\log p_{\bm{\alpha}}(\varsigma,N^a,N^b)\bigr]$ is the same term conditional on the price path $W=\varsigma$, and $b$ is a baseline independent of the trajectory. Since the score has zero mean, $\mathbb{E}_{\bm{\alpha}}[\nabla_{\bm{\alpha}}\log p_{\bm{\alpha}}]=0$, the baseline does not bias the estimator and only reduces its variance. In practice $b=\bar R_m$ is the batch mean of Algorithm~\ref{alg:reinforce_sigmm}.

\subsubsection{Sig-PPO (step-wise advantage).}
PPO \cite{schulman2017proximal} maximizes the same objective $J$, but factorises it through the per-step return decomposition $R_T=\sum_{j} r_{j+1}$, where $r_{j+1}$ is the reward accrued on $(\tau_j,\tau_{j+1}]$ (execution revenue and inventory penalty, with the terminal term charged at $\tau_M=T$), and replaces the path-wise score by per-step advantage estimates. Two ingredients worth discussing. First, as in Sig-Reinforce, the actor is learning and deploying a policy linear in the signature: this way we can separate the effect of how we are learning from what we are learning. Second, different from the other algorithm, a learned critic $V_\theta(\mathbb{Y}_{\tau_j})$, a neural network with parameters $\theta$, approximates the state value and yields the generalized advantage estimate
\[
    \widehat{A}_{j}
    \;=\;
    \sum_{h\ge 0} (\gamma\lambda_{\text{GAE}})^{h}\,
    \bigl[r_{j+h+1} + \gamma V_\theta(\mathbb{Y}_{\tau_{j+h+1}})
          - V_\theta(\mathbb{Y}_{\tau_{j+h}})\bigr].
\]
The actor is updated by gradient ascent on the clipped surrogate
\[
    L^{\text{clip}}(\bm{\alpha})
    \;=\;
    \mathbb{E}\!\left[\,
        \min\!\bigl(\rho_{j}\,\widehat{A}_{j},\;
                    \mathrm{clip}(\rho_{j},1{-}\epsilon,1{+}\epsilon)\,
                    \widehat{A}_{j}\bigr)
    \right],
\]
where
\[
    \rho_{j} \;=\;
    \frac{\pi_{\bm{\alpha}}(u_j\,|\,\phi_j)}
         {\pi_{\bm{\alpha}_{\text{old}}}(u_j\,|\,\phi_j)},
\]
Each rollout produces $M$ per-step gradients (one per event), reused across several inner optimization epochs over minibatches, the clipping acting as a trust-region proxy that bounds the policy drift per update. 
In the experiments we use $\gamma=1$ (no discounting over the finite horizon), $\lambda_{\text{GAE}}=0.95$, $\epsilon=0.2$, critic weight $0.5$, entropy weight $0.01$, six inner epochs with minibatches of size $\min(2048,\max(256,8B))$, gradient-norm clipping at $0.5$, Adam with learning rate $3\cdot10^{-4}$, and an MLP critic with $128$ hidden units; the critic additionally observes the price, which the actor does not.

\subsubsection{Key contrasts REINFORCE \textit{vs} PPO}
\begin{itemize}
\item \emph{Credit assignment. Sig-REINFORCE weights every score contribution $\nabla_{\bm{\alpha}}\log p(\text{event}_j)$, on both sides, by the same centered terminal reward $R_T-\bar R$; Sig-PPO weights each contribution by a step-local advantage $\widehat{A}_j$ estimated through the critic.}

\item \emph{Use of the structure. Sig-REINFORCE exploits the known, differentiable spread-to-intensity map through the explicit score \eqref{eq:poisson_score_alpha}; Sig-PPO is model-agnostic and treats the environment as a black box.}
\item \emph{Computational complexity.} Sig-REINFORCE evaluates a closed-form linear readout at each event and performs one closed-form update per batch; Sig-PPO must evaluate (and backpropagate through) a stochastic actor and a critic at every event, and then re-process the stored transitions over several epochs. On matched budgets this is an order of magnitude in wall-clock time (Section~\ref{sec:Poisson}).
\end{itemize}

\paragraph{Training configuration.}
Sig-REINFORCE uses $B=2048$ paths per rollout, $10^4$ iterations, and learning rates $5\cdot10^{-4}$ (Hawkes) and $10^{-4}$ (Poisson). Sig-PPO uses $B=2048$, $10^4$ iterations and the configuration described above. Both are evaluated out of sample on $16384$ fresh paths shared with the expert.

\subsection{Poisson order flow}\label{sec:Poisson}
As a controlled benchmark of our method, we adopt here the classical setting of \cite{avellaneda2008high,Dealingrisk}, Poisson arrivals with Bachelier price as described in subsection \ref{sec:numerics_setup}. Under a mean-variance criterion, the optimal market-making policy admits a known closed-form solution, see \cite[Chapter 6]{cartea2015algorithmic}. The question we aim to answer is whether our algorithm recovers the optimum in a regime where the analytical solution is available. Therefore, we assume that $N^a$ and $N^b$ are conditionally Poisson processes whose arrival intensities depend exponentially on the quoted spreads:
\[
\lambda_t^{a/b}(\bm{\alpha}^{a/b}) \;=\; A\,e^{-k\delta_{t^-}^{a/b}(\bm{\alpha}^{a/b})},\qquad t\in [0,T],
\]
for $\delta_{t^-}^{a/b}(\bm{\alpha}^{a/b}):=\langle\bm{\alpha}^{a/b},\mathbb{Y}_{t^-}\rangle$ and $A,k>0$.

Let $\mathbb{P}^0$ denote the reference measure under which both $N^a$ and $N^b$ are
unit-rate Poisson and independent of $W$. Let $\mathbb{F}=(\mathcal{F}_t)_{t\in[0,T]}$ be the natural augmented filtration of
$(W,N^a,N^b)$. By the Girsanov-Watanabe theorem for marked
point processes \cite[Thm.~7.4.I]{daley2003introduction}, the Radon Nikodym derivative reads as
\[
   \frac{d\mathbb{P}_{\bm{\alpha}}}{d\mathbb{P}^0}\bigg|_{\mathcal{F}_T}
   \;=\;
   \prod_{s\in\{a,b\}}
   \exp\!\left\{
       \int_0^T \log\lambda_t^s(\bm{\alpha}^s)\,dN_t^s
       - \int_0^T \bigl(\lambda_t^s(\bm{\alpha}^s) - 1\bigr)\,dt
   \right\}.
\]
Conditioning on $\mathcal{F}^W_T$, for any $\varsigma\in \mathcal C([0,T];\mathbb R)$,
\begin{align}\label{eq:loglik_factored}
\log p_{\bm{\alpha}}\bigl(\varsigma, N^a,N^b\,\big|\,\mathcal{F}^W_T\bigr)
= \int_0^T \log\lambda_t^a(\bm{\alpha}^a)\,dN_t^a - \int_0^T \lambda_t^a(\bm{\alpha}^a)\,dt
\nonumber\\
+\int_0^T \log\lambda_t^b(\bm{\alpha}^b)\,dN_t^b - \int_0^T \lambda_t^b(\bm{\alpha}^b)\,dt
+\mathrm{const}.\nonumber
\end{align}
Since
\[
   \nabla_{\bm{\alpha}}\log\lambda_t^{a/b}(\bm{\alpha}^{a/b})
   \;=\; -k\,\nabla_{\bm{\alpha}}\delta_{t^-}^{a/b}(\bm{\alpha}^{a/b}),
\]
\[
   \nabla_{\bm{\alpha}}\lambda_t^{a/b}(\bm{\alpha}^{a/b})
   \;=\; -k\,\lambda_t^{a/b}(\bm{\alpha}^{a/b})\,\nabla_{\bm{\alpha}}\delta_{t^-}^{a/b}(\bm{\alpha}^{a/b}),
\]
and, being  $\delta_{t^-}^{a/b}(\bm{\alpha}^{a/b})=\langle\bm{\alpha}^{a/b},\mathbb{Y}_{t^-}\rangle$, 
\[
   \nabla_{\bm{\alpha}^a}\delta_{t^-}^a \;=\; \nabla_{\bm{\alpha}^b}\delta_{t^-}^b \;=\;\pi_{\leq K}(\mathbb{Y}_{t^-}),\qquad
   \nabla_{\bm{\alpha}^a}\delta_{t^-}^b \;=\;
   \nabla_{\bm{\alpha}^b}\delta_{t^-}^a \;=\; 0,
\]
we finally get 
\begin{equation}\label{eq:poisson_score_alpha}
   \nabla_{\bm{\alpha}^{a/ b}}\log p_{\bm{\alpha}}
   = -k\!\int_0^T \!\pi_{\leq K}(\mathbb{Y}_{t^-})\,\bigl(dN_t^{a/b}-\lambda_t^{a/b}(\bm{\alpha}^{a/b})dt\bigr).
\end{equation}

\begin{remark}
Indexed by the event sequence $0=\tau_0<\tau_1<\cdots<\tau_M=T$ of
ordered ask and bid arrival times (with $\tau_M:=T$ as a terminal
sentinel), with inter-event intervals
\[
    \Delta\tau_{j} := \tau_{j+1}-\tau_{j},
    \qquad
    \Delta N_{j+1}^{a/b}\in\{0,1\}
\]
the indicator that the event at $\tau_{j+1}$ is of type $a$ or $b$, and
event-time features $\phi_{j}=\pi_{\leq K}(\mathbb{Y}_{\tau_{j}})$,
\eqref{eq:poisson_score_alpha} becomes
\[
    \nabla_{\bm{\alpha}^{a/b}}\log p_{\bm{\alpha}}
    \approx
    -k\sum_{j=0}^{M-1}\phi_{j}\big(\Delta N_{j+1}^{a/b}-\lambda_{j}^{a/b}\,\Delta\tau_{j}\big),
\]
where $\lambda_{j}^{a/b}:=\lambda^{a/b}(\delta_{\tau_{j}}^{a/b})$. Note that if $\lambda^a$ and $\lambda^b$ were
piecewise constant on each random inter-event interval
$[\tau_{j},\tau_{j+1})$, then $\lambda_{j}^{a/b}\Delta\tau_{j}$ would be the exact
compensator increment $\int_{\tau_{j}}^{\tau_{j+1}}\!\lambda_s^{a/b}\,ds$.
\end{remark}
We now compare numerically in this setting the Sig-REINFORCE algorithm with the Sig-PPO method described in Section \ref{sec:ppo}, and both of them against the analytical optimal solution of \cite{cartea2015algorithmic}. We instantiate the benchmark with the following values: horizon $T=300\,\mathrm{s}$, spread-coupling $k=0.2$, baseline order-arrival rate $A=0.517\,\mathrm{s}^{-1}$, price volatility $\sigma=2$ (in dollars per $\sqrt{\mathrm{s}}$) so that the Bachelier mid is $S_t=S_0+\sigma W_t$ with $S_0=100$ and a tick size of $0.01$. The risk-aversion parameter of the mean--variance reward is $\eta = 0.05$ and the per-side inventory cap is $Q_{\max}=10$. 

\begin{figure}[ht]
    \centering
    \includegraphics[width=\columnwidth]{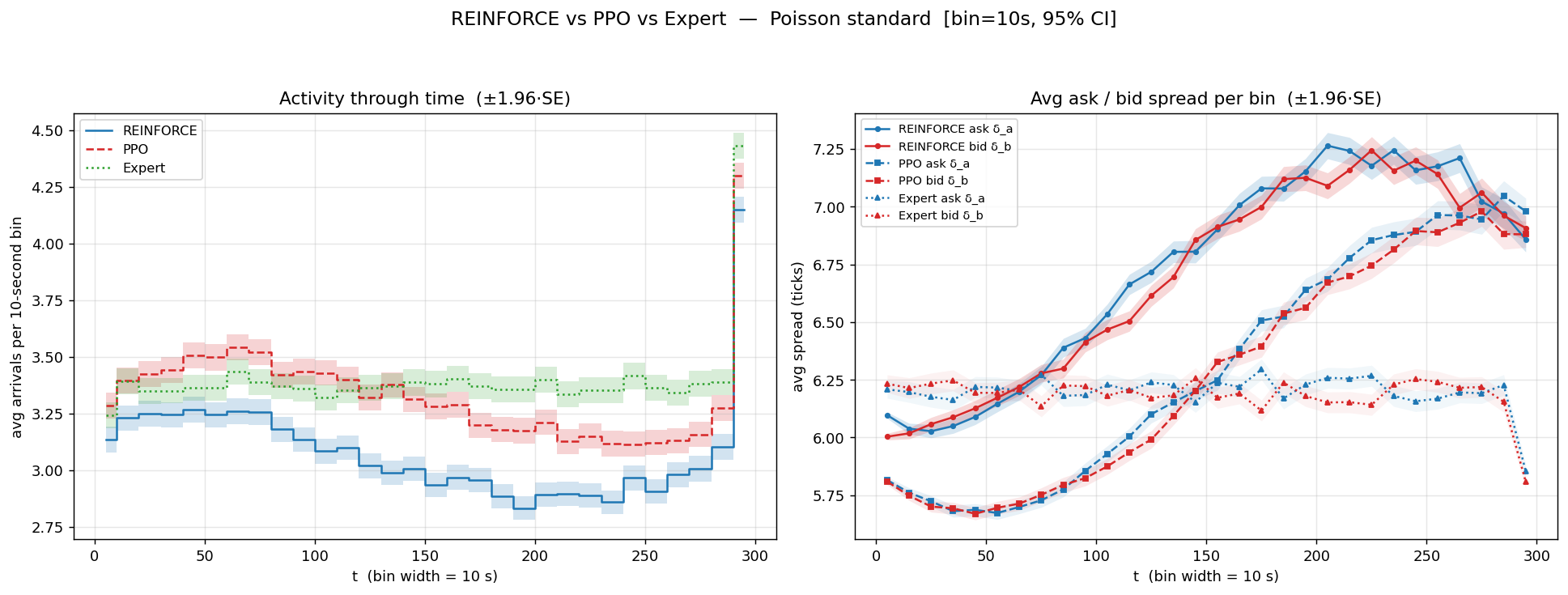}
    \caption{Averaged market order arrivals (left) and average bid and ask spread quoted (right) with Sig-REINFORCE, Sig-PPO and Expert (analytical solutions) for Poisson order flows.}
    \label{fig:poisson:flows}
    \Description{Two panels for Poisson order flow: average market order arrivals over time (left) and average bid and ask spreads quoted over time (right), comparing Sig-REINFORCE, Sig-PPO and the analytical expert.}
\end{figure}

\begin{figure}[ht]
    \centering
    \includegraphics[width=0.9\columnwidth]{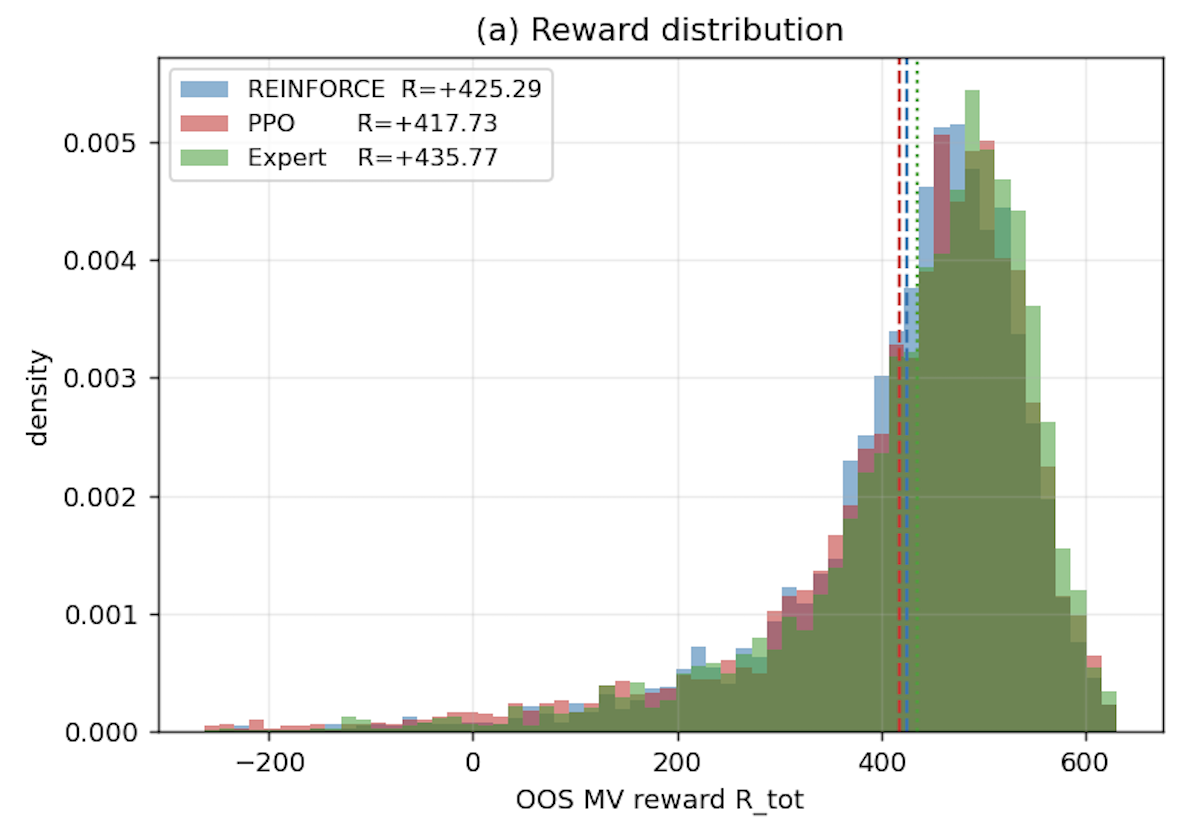}
    \caption{Payoff distribution for the market maker with Sig-REINFORCE, Sig-PPO and Expert for Poisson order flows.}
    \label{fig:poisson:value}
    \Description{Distribution of the mean-variance payoff for Sig-REINFORCE, Sig-PPO and the analytical expert under Poisson order flow, with the expert achieving the highest mean.}
\end{figure}

Figure \ref{fig:poisson:flows} compares, on the left, the average executed orders (total order flow) of the market maker and, on the right, their optimal quotes, obtained using Sig-REINFORCE, Sig-PPO, and the explicit analytical solution (the expert). We first observe a very similar pattern between Sig-REINFORCE and Sig-PPO relative to the benchmark. While the expert remains fairly stable, both Sig-REINFORCE and Sig-PPO exhibit more variance in the best ask/bid quotes. Notably, Sig-REINFORCE achieves this with a slightly lower quote variance than Sig-PPO. The bid/ask spreads are also wider under Sig-REINFORCE and Sig-PPO, consequently leading to a lower average total order flow. 
We now turn to the reward (payoff) of the mean-variance problem. We simulate it in Figure \ref{fig:poisson:value} for Sig-REINFORCE, Sig-PPO, and the expert, and we also report the value of the problem at time $300$ in the box at the top right. While the reward obtained by the expert is, as expected, the highest, the gap with Sig-REINFORCE and Sig-PPO remains small: over common evaluation paths the mean rewards are $425.3\pm2.1$ for Sig-REINFORCE and $417.7\pm2.5$ for Sig-PPO (standard errors), against $435.8$ for the expert. The difference between the two learning methods is about two standard errors, which we read as parity with a slight edge for Sig-REINFORCE rather than as a systematic advantage; the practical difference between the two is computational. On the matched budget ($2048$ paths per iteration, $10^4$ iterations) training Sig-REINFORCE took (on average of $10$ experiments) $4100\,$s against more than $30000\,$s for Sig-PPO. The gap stems from evaluating a stochastic actor and a critic network at every event of every path, while Sig-Reinforce does computations once per path at the end. Timings were obtained on a single MacBook Air (Apple M-, 4P+4E cores, 16 GB), CPU only.

\subsection{Hawkes order flow}\label{sec:hawkes_numerics}
We now turn to the case where the order flow follows a Hawkes process. Preliminary works on market making with self-exciting order flows have been investigated in \cite{jusselin2021optimal}, but the theoretical study of the mean-variance problem remains open. This section focuses on the numerical analysis of this problem with Sig-REINFORCE algorithm compared with Sig-PPO. We use the Br\'emaud-Massouli\'e representation (thinning) \cite{ogata1981lewis,bremaud1981point} to represent the Hawkes processes.
Let $\Pi^{a},\Pi^b$ be two Poisson random measures on $\mathbb R^+\times\mathbb R^+$ with compensators $\tilde \Pi^{a/b}(ds,d\theta)=\Pi^{a/b}(ds,d\theta)-ds\,d\theta$. We define
\[N^{a/b}_t:=\int_0^t\int_0^\infty\mathbf 1_{\theta\leq \lambda^{a/b}_s}\Pi^{a/b}(ds,d\theta),\] with intensity $\lambda_t^{a/b}:=\lambda^{a/b}(t,\delta^{a/b}_{t^-})$ given by 
\[
\lambda^{a/b}(t,\delta^{a/b}_{t^-}) = e^{-k \delta_{t^-}^{a/b}}(\mu_{a/b}+X^{a/b}_t),
\]
where,
\[
X^{a/b}_t = \int_0^te^{-\beta (t-s)} \int_0^\infty   \mathbf 1_{\theta\leq \lambda^{a/b}_s} \Pi^{a/b}(ds,d\theta). 
\] 

In the numerical experiments we choose $T=300$, $k=0.2$ and $\mu_a=\mu_b=0.3102$, $\beta=2$.

\begin{remark}
The likelihood-ratio argument extends naturally to the Hawkes setting through
the Br\'emaud--Massouli\'e representation, which
constructs the controlled Hawkes process as a thinning of a Poisson random
measure with predictable intensity
\[
\lambda_t^{a/b}
=
e^{-k\delta_{t^-}^{a/b}}
\bigl(\mu_{a/b}+X_t^{a/b}\bigr).
\]
The Br\'emaud-Massouli\'e representation is essential since it realizes the
controlled Hawkes process and its intensity as the strong solutions of a system of stochastic equations
driven by a Poisson random measure (independent of $\bm{\alpha}$). In particular, the intensity processes
$\lambda^{a/b}$
are predictable with respect to the filtration generated by the underlying
Poisson random measures $\Pi^{a/b}$. Consequently, for every realized trajectory
$\omega$, the intensity is an adapted functional of the observed history,
so that the standard likelihood of a point process with predictable
intensity applies pathwise. For a realized trajectory
$\omega=(\tau_1,\tau_2,\ldots,\tau_M)$,
the jump times are fixed, so that
$
X_t^{a/b}(\omega)
=
\sum_{\tau_{j}<t}
e^{-\beta(t-\tau_{j})}
$
is a predictable process determined by the observed history.
Moreover, by definitions of the arrival times
$\tau_{j},\tau_{j+1}$, no new jump occurs in the interval
$(\tau_{j}, \tau_{j+1})$ and the excitation evolves only through exponential decay,
namely
\[
X_t^{a/b}
=
e^{-\beta(t-\tau_{j})}
X_{\tau_{j}}^{a/b},
\qquad
t\in[\tau_{j},\tau_{j+1}).
\]
Consequently, the compensator over this interval is
\[
\Delta\Lambda_{j}^{a/b}
=
\int_{\tau_{j}}^{\tau_{j+1}}
\lambda_t^{a/b}\,dt
\approx
e^{-k\delta_{j}^{a/b}}
\left(
\mu_{a/b}\Delta\tau_{j}
+
X_{\tau_{j}}^{a/b}
\frac{1-e^{-\beta\Delta\tau_{j}}}{\beta}
\right).
\]

Sig-REINFORCE estimates the policy gradient when the score is computed for a fixed realized trajectory
$\omega$. Once the trajectory has been sampled, the jump times and hence
$X_{\tau_{j}}^{a/b}(\omega)$ are observed data, while the parameter
dependence enters only through the control
$\delta_{j}^{a/b}=\langle \bm{\alpha}^{a/b}, \mathbb{Y}_{\tau_j} \rangle$.
Therefore two separate identities hold. First, since $\delta_j^{a/b}=\langle\bm{\alpha}^{a/b},\mathbb{Y}_{\tau_j}\rangle$ and the excitation $X^{a/b}_{\tau_j}$ is a functional of the observed jump times only, setting $\phi_{j}=\pi_{\leq K}(\mathbb{Y}_{\tau_{j}})$,
\[
\nabla_{\bm{\alpha}^{a/b}}\log\lambda_j^{a/b} = -k\,\phi_j .
\]
Second, because $\delta_j^{a/b}$ is held fixed on $[\tau_j,\tau_{j+1})$ the compensator increment is proportional to
$e^{-k\delta_j^{a/b}}$, hence
\[
\nabla_{\bm{\alpha}^{a/b}}\Delta\Lambda_j^{a/b}
 = -k\,\phi_j\,\Delta\Lambda_j^{a/b},
\]
exact for the piecewise-constant quotes of the discretisation and first-order accurate for the true compensator. Combining the two in the log-likelihood gives
\[
 \nabla_{\bm{\alpha}^{a/b}}\log p_{\bm{\alpha}}
= -k\sum_{j=0}^{M-1}
\phi_{j}
\left(
\Delta N_{j+1}^{a/b}
-
\Delta\Lambda_{j}^{a/b}
\right).
\]
\end{remark}

\begin{figure}[htbp]
        \centering
        \includegraphics[width=0.9\columnwidth]{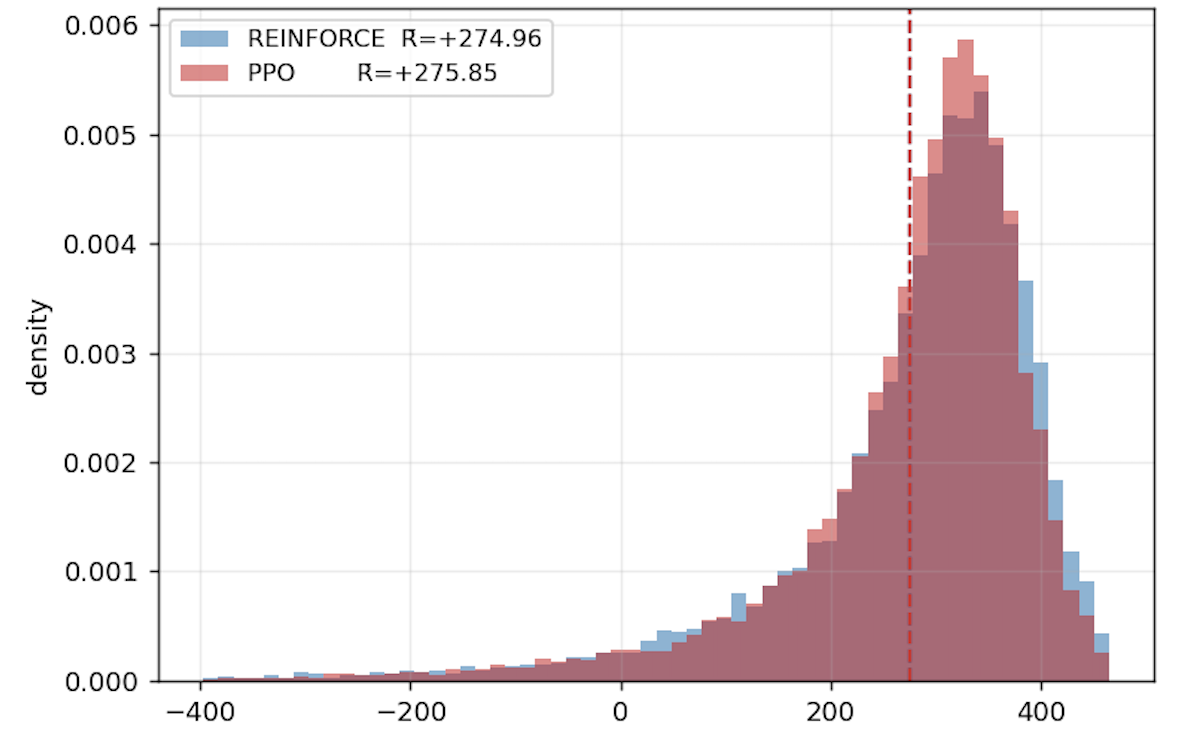}
        \caption{Return Sig-REINFORCE vs. Sig-PPO}\label{hawkes:reward}
    \Description{Distribution of the mean-variance return under Hawkes order flow for Sig-REINFORCE and Sig-PPO, which are statistically close.}
\end{figure}
\begin{figure}[htbp]
        \centering
        \includegraphics[width=\columnwidth]{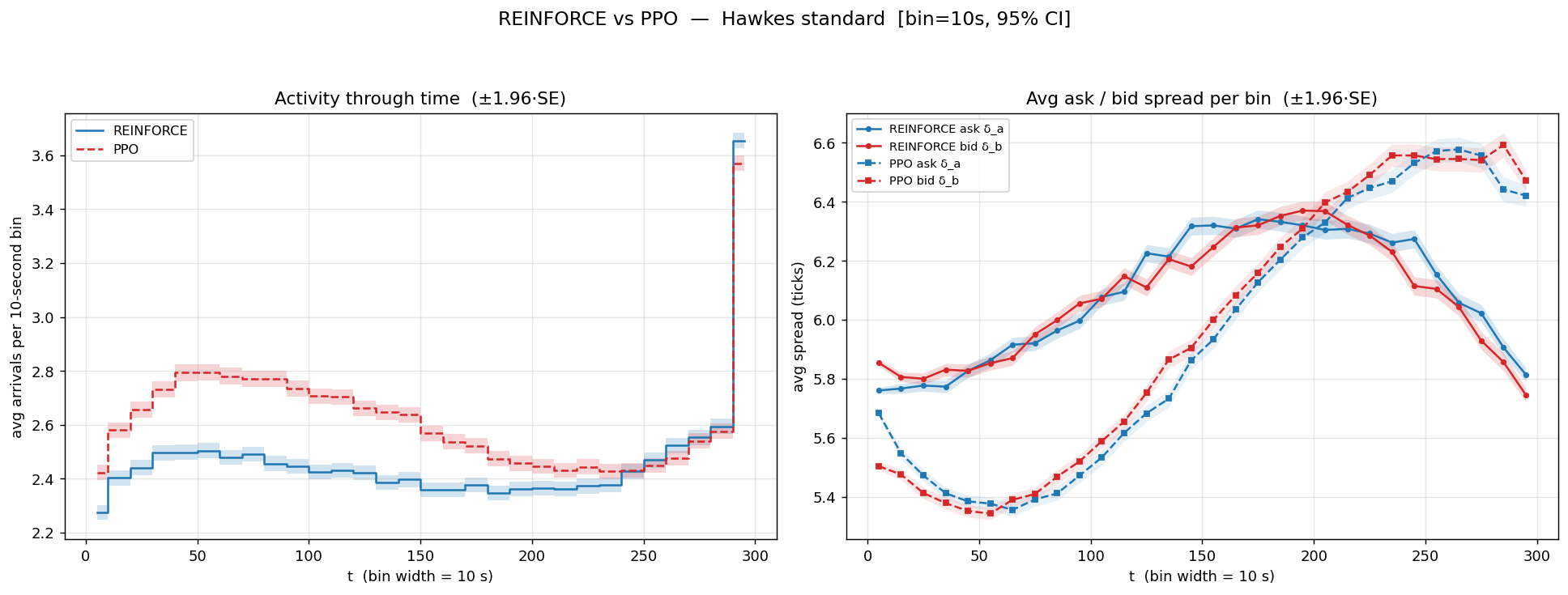}
        \includegraphics[width=\columnwidth]{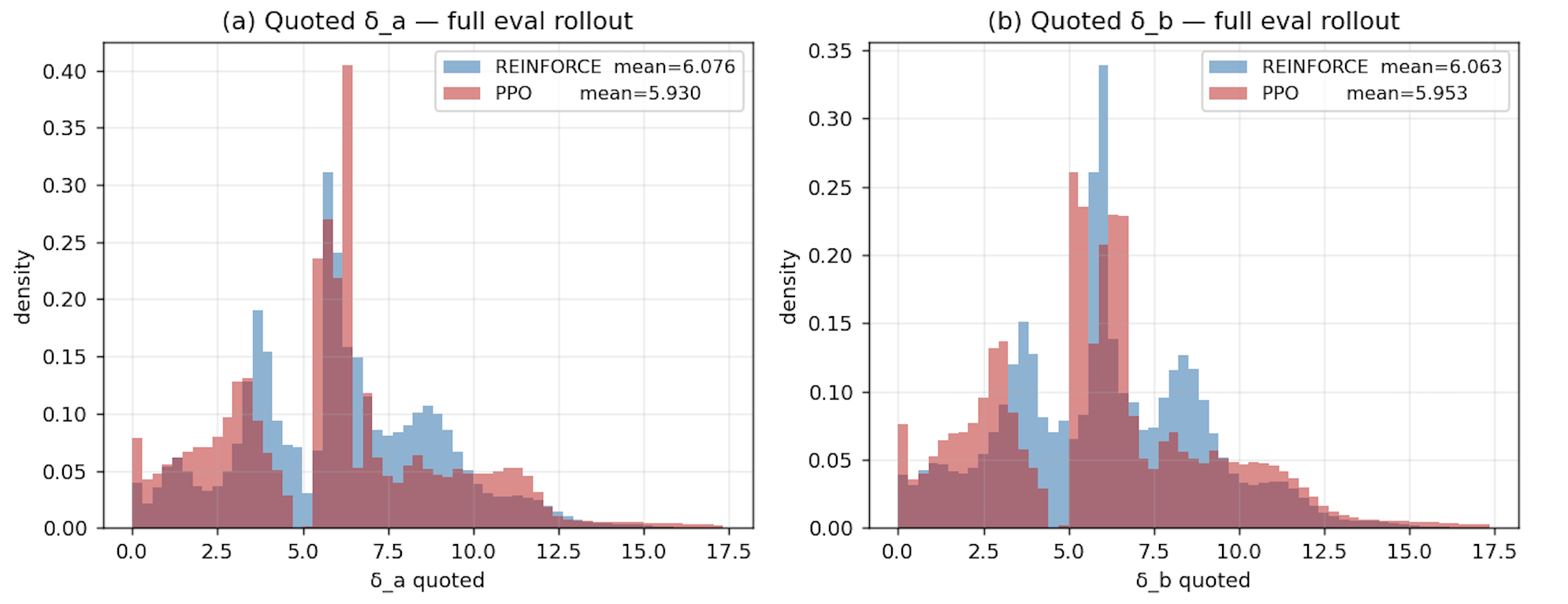}
        \caption{Expected Bid-Ask spread value  with Sig-REINFORCE vs. Sig-PPO (top) and respective distributions (bottom)}\label{Hawkes:bidask average}
    \Description{Top: expected bid-ask spread over time for Sig-REINFORCE and Sig-PPO under Hawkes order flow. Bottom: distributions of the quoted spreads for the two methods.}
\end{figure}
The reward distribution associated to the mean-variance problem are presented in Figure \ref{hawkes:reward}. The average reward with Sig-REINFORCE is statistically indistinguishable from Sig-PPO, which shows how it is possible to learn equally well looking at entire paths. Digging deeper into the policy found by the two algorithms, Figure \ref{Hawkes:bidask average} compared the bid-ask spread obtained with our Sig-REINFORCE method and the one derived from Sig-PPO. We first observed that on a reasonable trading scale of $T=300s$, Sig-REINFORCE has slightly less variance that Sig-PPO, especially close to the end of the horizon, similarly to the Poisson case previously investigated. This is confirmed by the distribution graphs which seems to be a bit more bimodal for Sig-PPO compared with Sig-REINFORCE. The average spread is again slightly higher with Sig-REINFORCE compared with Sig-PPO, but the shapes of the average quote in time for both algorithm is consistent.

\section{Conclusion}
This paper develops a signature-based framework for optimal market making with controlled order arrivals and introduces the Sig-REINFORCE algorithm. In the Poisson setting, it achieves a reward close to the analytical benchmark; in the Hawkes setting, it matches Sig-PPO with lower variance and at a much lower computational cost. These results provide a proof of concept for signature methods in stochastic control with jumps.

Several limitations remain. We provide no convergence guarantees or error bounds in terms of the truncation level, and the interplay between truncation, estimation and optimization errors is open, for this problem and more generally for jump diffusions (see, e.g., \cite{guo2025consistency,bank2025stochastic}). Higher truncation levels increase expressiveness but also computational and statistical cost, and controlling the variance of the gradient estimator requires care even for a single asset. Better variance reduction and feature selection would allow extensions to multiple assets, correlated prices and mutually exciting order flows.

A natural extension of the price dynamics is rough volatility \cite{gatheral2022volatility}. Although the price can remain a semimartingale, adding fractional Brownian or Volterra volatility factors to the signature state would require joint signature representations with the jump order flow, together with integrability properties needed for the control problem.

\begin{acks}
We thank the anonymous referees for their insightful comments and suggestions, which helped improve the paper.
\end{acks}

\bibliographystyle{ACM-Reference-Format}
\bibliography{myref}

\end{document}